\documentclass[a4paper,leqno]{article}

\usepackage[T1]{fontenc}
\usepackage[latin1]{inputenc}                                                                     
\usepackage{mathtools,amssymb,amsthm}
\mathtoolsset{mathic}
\usepackage{lmodern}
\usepackage{eucal}
\usepackage{hyperref}
\usepackage{microtype}
\usepackage{cleveref}

\newlength{\alphabet}
\usepackage[textwidth=3\alphabet]{geometry}

\theoremstyle{plain}
\newtheorem{thm}{Theorem}[subsection]
\newtheorem{prop}[thm]{Proposition}
\newtheorem{coro}[thm]{Corollary}
\newtheorem{lem}[thm]{Lemma}

\newtheorem{Conj}[thm]{Fried' s Conjecture}

\newtheorem{defi}[thm]{Definition}
\newtheorem{rmrk}[thm]{Remark}
\newtheorem*{ack}{Acknowledgements}

\newcommand{\R}{\mathbb{R}}
\newcommand{\Z}{\mathbb{Z}}

\renewcommand{\H}{\mathbb{H}}
\newcommand{\C}{\mathbb{C}}

\newcommand*{\rom}[1]{\romannumeral}
\renewcommand{\Re}{{Re}}

\newcommand{\geod}{\mathrm{geod}}

\DeclareMathOperator{\PSL}{PSL}
\DeclareMathOperator{\GL}{GL}

\DeclareMathOperator{\PSO}{PSO}
\DeclareMathOperator{\SU}{SU}

\DeclareMathOperator{\Spin}{Spin}
\DeclareMathOperator{\Sp}{Sp}
\DeclareMathOperator{\F}{F}

\DeclareMathOperator{\RE}{Re}
\DeclareMathOperator{\Tr}{Tr}
\DeclareMathOperator{\tr}{tr}
\DeclareMathOperator{\End}{End}
\DeclareMathOperator{\Vol}{Vol}

\DeclareMathOperator{\Id}{Id}
\DeclareMathOperator{\Ad}{Ad}

\DeclareMathOperator{\Norm}{Norm}

\DeclareMathOperator{\spec}{spec}

\DeclareMathOperator{\sign}{sign}

\DeclareMathOperator{\prim}{prime}

\DeclareMathOperator{\gr}{gr}
\DeclareMathOperator{\Rep}{Rep}

\DeclareMathOperator{\Ker}{Ker}

\DeclareMathOperator{\tor}{tor}
\DeclareMathOperator{\mult}{mult}
\DeclareMathOperator{\Fix}{Fix}

\numberwithin{equation}{section}
\numberwithin{thm}{section}

\renewcommand{\thethm}{
\ifnum\value{subsection}=0
\thesection
\else
\thesubsection
	\fi
		.\arabic{thm}
}

\title{Twisted dynamical zeta functions and the Fried's conjecture}
\author{Polyxeni Spilioti}

\begin{document}
\maketitle

\abstract {This is a survey article on the twisted dynamical zeta functions of Ruelle and Selberg
 and the Fried's conjecture. It is based on the mini-course: "Twisted Ruelle zeta function, complex-valued analytic torsion and the Fried's conjecture", 
given by the author during the thematic trimester programme:
"Representation Theory and Noncommutative Geometry"
at the Institut Henri Poincar\'{e}.}

\tableofcontents

\section{Introduction}

\subsection{Dynamical zeta functions}

Initially, we consider the dynamical zeta functions of \textit{Selberg and Ruelle} for compact hyperbolic surfaces $X$. These dynamical-geometrical objects are defined in terms of the lengths of the \textit{closed geodesics} on $X$.  A \textit{prime closed geodesic} is a closed geodesic that can not be written as the multiple of a shorter closed geodesic.

According to Ruelle \cite{ruelle2002dynamical}, the grandmother of all zeta functions is the \textit{Riemann zeta function}
\begin{equation*}
\zeta(s)=\sum_{n=1}^{\infty} n^{-s},\quad \RE(s)>1.
\end{equation*}
Completely analogously to the Euler product formula for the Riemann zeta function, 
\begin{equation*}
\zeta(s)=\prod_{p\hspace*{0.6pt}\text{prime}}(1-p^{-s})^{-1}, \quad \RE(s)>1,
\end{equation*}
the Ruelle zeta function is given by an Euler product formula,
where now the product runs over the prime closed geodesics
of length $l(\gamma)$,
\begin{equation}\label{Ruelle1}
R(s)=\prod_{\gamma\text{ prime}}(1-e^{-sl(\gamma)}), \quad \RE(s)>1.
\end{equation}
Correspondingly, the Selberg zeta function is defined  by a double product
\begin{equation}\label{Selberg1}
Z(s) = \prod_{\gamma\text{ prime}}\prod_{k=0}^\infty\det\left(1-e^{-(s+k)\ell(\gamma)}\right),
\quad \RE(s)>1.
\end{equation}
One can view the Ruelle zeta function as the \textit{geometrical analogue} of the Riemann zeta function and the precise analogy is
\begin{equation*}
p \hspace*{0.9pt}\text{prime}\longleftrightarrow e^{l(\gamma)}.
\end{equation*}

Even more interesting is the fact that there are precise analogies between the Riemann zeta function and the Selberg zeta function, concerning their analytic properties. 
Namely, the Selberg zeta function has a meromorphic continuation to the whole complex plane and there exists the so called \textit{spectral zeros}, which have real part equal to $\frac{1}{2}$ (\cite{Selberg1956Harmonic}).
The key point to prove this result is the \textit{Selberg trace formula}, which was introduced by Selberg in \cite{Selberg1956Harmonic}. In this work, Selberg already pointed out the analogy of the trace formula to formulas arising in number theory from the \textit{zeta- and $L$-functions of algebraic number fields} as well as the analogue of the \textit{Riemann hypothesis} for the Selberg zeta function.

Apart form the above analogies, the Selberg zeta function provides a link to number theory via the \textit{trace formula}. In addition, the dynamical zeta functions on arithmetic hyperbolic surfaces or manifolds in general can be very interesting objects to study, see for example the work of Sarnak and Xue (\cite{SarnakXue}) and Marshall and M{\"u}ller (\cite{MMar}).
There are also several other zeta functions from number theory that could be considered as analogous to the dynamical zeta functions. Moreover, there are several analogies concerning formulas, results and and conjectures. We refer the reader to the very interesting article of 
Deninger, where these analogies are set in place for the \textit{Hasse-Weil zeta function} and the \textit{Lichtenbaum's conjecture} (\cite{deninger2007dynamicalsystemsanaloguelichtenbaums}).

On the other hand, the Ruelle zeta function associated with a dynamical system 
was introduced by Ruelle (\cite{ruelle2002dynamical}), a mathematical physicists and his definition is related to \textit{thermodynamic formalism}.
As it is mentioned before, the singularities of the Selberg zeta function, a dynamical zeta function associated with the geodesic flow on a hyperbolic surface, can be described in terms of the spectrum of the Laplace operator and hence there is a strong link to the spectral theory of the elliptic, differential operators. From mathematical physics point of view, the connection between the goedesic flow-classical mechanics and the Laplace operator-quantum mechanics is related to the so-called  \textit{quantum chaos} (\cite{ruelle2002dynamical}).

In other words, the dynamical zeta functions can be related to problems from \textit{geometry, topology, dynamics, spectral geometry, number theory and mathematical physics} and that is exactly what makes them a fascinating object to study.

\subsection{Fried's conjecture}

The first question one may ask is what are the analytic properties of the Ruelle zeta function.
By \eqref{Ruelle1} and\eqref{Selberg1}, one can easily see that
\begin{equation*}
R(s)=\frac{Z(s)}{Z(s+1)}
\end{equation*}
Are there special values for $R(s)$?
Secondly, one may want to consider representations of the lattice $\Gamma$, since 
\textit{holonomy representations}  are often considered.
Moreover, these dynamical zeta functions can be studied under different and more
complicated geometrical considerations, for example, higher-dimensional hyperbolic 
manifolds with cusps, orbifolds, general locally symmetric spaces and
Riemannian manifolds of variable negative curvature. 

\begin{Conj}\label{Fcon}
\textit{The Fried's conjecture is the conjecture that asks if the special value of the Ruelle zeta function at zero is related to a topological or spectral invariant.}
\end{Conj}

In \cite{fried1987lefschetz}, Fried posed,\textit{ among several other questions}, the question of finding such a relation, when one considers \textit{arbitrary representations} of the fundamental group of an oriented Riemannian manifold of constant negative curvature (\cite[p.62-63]{fried1987lefschetz}).
Fried, distinguished two cases, where the dimension of the manifold is $2$ or greater.

If the representation is orthogonal, both cases were proved by Fried in \cite{fried1986fuchsian} and
\cite{Fried}.
Let $X$ be a closed, oriented, hyperbolic manifold $X$ of dimension $d$.
Let $\rho\colon\pi_{1}(S(X))\rightarrow O(m)$ be an acyclic representation of 
$\pi_{1}(SX)$.

The twisted Ruelle zeta function $R_{\rho}(s)$, associated with $\rho$,
is defined for $\RE(s)>d-1$ by
\begin{equation*}
R_{\rho}(s)=\prod_{\gamma}\det(I-\rho(\gamma)e^{-sl(\gamma)}),
\end{equation*}
where  $\gamma$ runs over the prime closed geodesics of $X$, i.e. those not expressible as multiples of some shorter closed geodesic and $l(\gamma)$ denotes the length of $\gamma$.

\begin{thm}(Fried \cite[Theorem 1]{Fried})
The twisted Ruelle zeta function $R_{\rho}(s)$
extends meromorphichally to $\C$, and for $\epsilon=(-1)^{d-1}$,
\begin{equation*}
\vert  R_{\rho}(0)^{\epsilon}\vert=\tau_{\rho}(S(X)).
\end{equation*}
\end{thm}
Here, $\tau_{\rho}(S(X))$ denotes the Reidemeister torsion of $S(X)$ 
(\cite{DR,franz1935torsion,reidemeister1935homotopieringe}).

In the case of an \textit{arbitrary representation} of the fundamental group of a
compact hyperbolic manifold, the resolution of the Fried's conjecture \ref{Fcon}
was given by the author, under certain assumptions (\cite{spilioti2025twisted}) by Frahm and the author (\cite{frahm2023twisted}) and by Frahm, B{\'e}nard and the author(\cite{JEP_2023__10__1391_0}).
More precisely, we have the following theorems.

\subsubsection{Compact hyperbolic surface}
 
 	\begin{thm}(\cite[Theorems 4.2.6 and 4.2.9]{frahm2023twisted})
		Let $X=\Gamma\backslash \H^2$ be a compact hyperbolic surface and let $\chi\colon \Gamma\rightarrow \GL(V_{\chi})$ be a finite-dimensional, complex representation of $\Gamma$. Then, the twisted Selberg zeta function $Z(s;\chi)$ admits a meromorphic continuation to $\C$ which satisfies certain functional equations.
	\end{thm}

	\begin{coro}(\cite[Corollary 4.2.8]{frahm2023twisted})
		The twisted Ruelle zeta function $R(s;\chi)$ admits a meromorphic continuation to $\C$ which satisfies the functional equation
		$$ R(s;\chi)R(-s;\chi)=(2\sin\pi s)^{2(2g-2)\dim V_\chi}. $$
	\end{coro}
	
	\begin{coro}(\cite[Corollary 4.2.11]{frahm2023twisted})
		The behavior of the twisted Ruelle zeta function $R(s;\chi)$ near $s=0$	is given by
		\begin{equation*}
			R(s;\chi)=\pm(2\pi s)^{\dim(V_{\chi})(2g-2)}+\text{higher order terms}.
		\end{equation*}
	\end{coro}

\subsubsection{Compact hyperbolic orbisurface}\label{Friedodd}

Let $X$ be a compact hyperbolic surface $X$ with finitely many singular points $x_1, \ldots, x_r$ of finite order. Let $X_1$ be its unit tangent bundle. The generator $\mathfrak X$ of the geodesic flow acts on $X_1$ and its (prime) periodic orbits are lifts of (prime) closed geodesics on the surface.
For each singular point $x_j$ on $X$ there is a class of loops $c_j$ encircling $x_j$, which has finite order $\nu_j \in \Z_{>0}$ in the orbifold fundamental group of $X$. 
Their lifts (still denoted by $c_j$) to $\pi_1(X_1)$ satisfy $c_j^{\nu_j}= u$,  
where $u$ is the class of the generic fiber. The latter represents a loop in $X_1$ of a unit vector which makes a full positive rotation along the fiber around the base point.

Let $\rho \colon \pi_1(X_1) \to \GL(V_\rho)$ be a representation with $\dim V_\rho=n$. For any $j=1, \ldots, r$ we denote by $n_j= \dim \Fix  \rho(c_j)$, so that $\rho(c_j) = I_{n_j} \oplus T_j$.

\begin{thm} (\cite[Theorem A]{JEP_2023__10__1391_0})\label{orbispecialvalue}
	For any  irreducible representation $\rho\colon \pi_1(X_1) \to \GL(V_\rho)$, the Ruelle zeta function $R(s;\rho)$ converges on some right half plane in $\C$ and extends meromorphically to the whole complex plane. Moreover:
\begin{enumerate}
\item \label{i1}
If $\rho(u)= \Id_{V_\rho}$, then $R(s;\rho)$ vanishes at $s=0$ with order and leading coefficient prescribed by:
$$R\left(\frac{s}{2\pi}, \rho\right) \sim_{s \to 0}  \pm \frac{s^{n(2g-2+r)-\sum_{j=1}^rn_j}}{\prod_{j=1}^r  |\det(I_{n-n_j} -T_j)|\,(-\nu_j)^{-n_j}}.$$

\item \label{i2}
If $\rho(u) \neq\Id_{V_\rho}$, then the representation $\rho$ is acyclic. Let $\mathfrak e_{\geod}$ be the Euler structure induced by the geodesic flow on $X_1$. Then
$$R(0;\rho) = \pm\tor(X_1,\rho,\mathfrak {e}_{\geod}, \omega^1)$$
 where $\tor(X_1,\rho, \mathfrak e_{\geod}, \omega^1)\in\C^\times$ denotes the Reidemeister--Turaev torsion of $X_1$ in the representation $V_\rho$, the Euler structure $\mathfrak e_{\geod}$ and the natural homology orientation $\omega^1$ (see \cite[Section 2]{JEP_2023__10__1391_0}
for details).
\end{enumerate}
The sign is in both cases equal to $(-1)^{\mult(0;\Delta_{\tau,\rho}^\sharp)}$, where $\mult(0;\Delta_{\tau,\rho}^\sharp)$ is the multiplicity of the eigenvalue $0$ of the twisted Bochner--Laplacian $\Delta_{\tau,\rho}^\sharp$ acting on sections of a certain orbifold bundle over $X_1$
(see \cite[Section 3.1 and 3.2]{JEP_2023__10__1391_0} for details).
\end{thm}

\subsubsection{Odd dimensional compact hyperbolic manifolds}

Let $X=\Gamma\backslash \H^{d}$ be compact real hyperbolic manifold of odd dimension $d$.
Let $\chi\colon\Gamma\rightarrow \GL(V_{\chi})$ be a finite-dimensional representation of $\Gamma$. Let $\widehat{M}$ be the the set of equivalence classes 
of unitary and irreducible representations of the group $M$
(see Sections \ref{lss1} and \ref{tdzf}).
We consider the twisted Ruelle zeta function $R(s;\sigma,\chi)$,
defined as in Definition \ref{ruelledefiodd}.

We have the following theorem.

\begin{thm} (\cite[Theorem 7.13]{Spilioti2018})
For every $\sigma\in\widehat{M}$, the Ruelle zeta function $R(s;\sigma,\chi)$ admits a meromorphic
continuation to the whole complex plane $\C$.
\end{thm}

In \cite{spilioti2020functional}, a \textit{determinant formula} is derived, 
which is a representation of the twisted Ruelle zeta function as a 
\textit{product of regularized determinants of  certain elliptic operators}, induced by the \textit{twisted Laplacians}. By this formula, we have the following
theorem.

\begin{thm} (\cite[Theorem 1]{spilioti2025twisted})
Let $\chi$ be a finite-dimensional complex representation of $\Gamma$.
Let  $\Delta^{\sharp}_{\chi,k}$ be the flat Hodge Laplacian, acting on 
the space of $k$-differential forms on $X$ with values in the flat vector bundle $E_{\chi}$.
Then, the Ruelle zeta function has the representation
\begin{align}\label{detmain}
R(s;\chi)=\notag&\prod_{k=0}^{d-1}\prod_{p=k}^{d-1}{{\det}_{\gr}\big(\Delta^{\sharp}_{\chi,k}+s(s+2(\vert\rho\vert-p))\big)}^{(-1)^{p}}\\
&\cdot\exp\bigg((-1)^{\frac{d-1}{2}+1}\pi(d+1)\dim (V_{\chi})\frac{\Vol(X)}{\Vol(S^{d})}s\bigg),
\end{align}
where $\Vol(S^{d})$ denotes the volume of the $d$-dimensional Euclidean unit sphere.
Let $d_{\chi,k}:=\dim\Ker(\Delta_{\chi,k})$. Then, the singularity of the Ruelle zeta function at $s=0$ is of order
\begin{equation*}
\sum_{k=0}^{(d-1)/2}(d+1-2k)(-1)^{k}d_{\chi,k}.
\end{equation*}
\end{thm}

One can assume that the representation is acyclic, and hence the cohomologies vanish,
but \textit{Hodge theory is no longer applicable}: there is no isomorphism between the cohomologies
with coefficients in the local system defined by $\chi$,
$H^{k}(X;E_{\chi})$,
and the kernels of the corresponding non-self-adjoint 
twisted Hodge Laplacians,
$\Ker(\Delta^{\sharp}_{\chi,k})$, i.e., 
\begin{equation*}
 \Ker(\Delta^{\sharp}_{\chi,k})\ncong H^{k}(X;E_{\chi}).
\end{equation*}
Hence, the kernels might \textit{not} be trivial
and we can not conclude by the determinant formula \eqref{detmain}
that the twisted Ruelle zeta function is regular at zero.

To overcome this problem, we considered in \cite{spilioti2025twisted}
acyclic representations of $\Gamma$, which are non-unitary
but close enough to unitary ones.
By a \textit{continuity argument} (see \cite[Proposition 6.8]{BK2}),
the twisted Hodge Laplacian is injective and 
by the determinant formula
one can conclude the regularity of the twisted Ruelle zeta function at zero. 
Again by the determinant formula, we found an \textit{explicit relation} 
between the twisted Ruelle zeta function at zero and the 
\textit{complex-valued analytic torsion} (defined as in \cite{cm}) and the 
\textit{refined analytic torsion} (defined as in \cite{BK2}).

Let $\Rep(\pi_{1}(X),\C^{n})$  be the set of all $n$-dimensional, complex representations of $\pi_{1}(X)$. This set has a natural structure of a \textit{complex algebraic variety}.
Let  $V \subset \Rep(\pi_{1}(X),\C^{n})$ be an 
\textit{open neighbourhood} (in classical topology) of the set  $\Rep^{u}_{0}(\pi_{1}(X),\C^{n})$
of acyclic, unitary representations such that, for all $\chi\in 
V$, the \textit{odd signature operator} $B_{\chi}$ is bijective
(see \cite[Section 3]{spilioti2025twisted} and \cite{BK2}).

We have the following theorem.

\begin{thm}(\cite[Theorem 2]{spilioti2025twisted})
Let $\chi\in V$. Then, the Ruelle zeta function $R(s;\chi)$ is regular at $s=0$ and is equal 
to the complex Cappell-Miller torsion, 
\begin{equation*}
R(0;\chi)=\tau_{\chi}.
\end{equation*}
\end{thm}

\subsubsection{Related work}

The last years, there has been much and spectacular progress in resolution of this conjecture by several mathematicians under different settings and 
using different techniques. 

Bunke-Olbrich  in \cite{BO} studied the twisted dynamical zeta functions 
associated with unitary representations for all locally symmetric spaces 
of real rank one and proved  the Fried's conjecture.
In the higher rank case and for unitary representations the Fried's conjecture was proved by 
Shen (\cite{shen2017analytic}), see also Moscovici-Stanton (\cite{MS91}),
Moscovici-Stanton-Frahm
(\cite{MSF}).
For other contributions to Fried's conjecture
and also applications of its resolution, see M\"{u}ller
(\cite{muller2012asymptotics}, \cite{MMar}, \cite{muller2021ruelle}),
Pfaff (\cite{pfaff2012selberg}),
Shen (\cite{shen2021analytica}, \cite{shen2021analytic}, \cite{shen2023complex}), 
Shen-Yu (\cite{MR4530056}),
Bismut-Shen (\cite{bismut2025anosov}),
Deitmar (\cite{deitmar2011holomorphic}),
Hochs-Saratchandran
(\cite{hochs2025ruelledynamicalzetafunction}),
Hochs-Pirie (\cite{hochs2025equivariantfriedconjecturesuspension}),
Dyatlov-Zworski (\cite{dyatlov2017ruelle}),
Guillarmou-Dang-Rivi\`{e}re, Shen (\cite{dang2020fried}),
Chaubet-Dang (\cite{chaubet2024dynamical}) and Ceki\'{c}-Delarue-Dyatlov-Paternain (\cite{cekic2022ruelle}), Jorgenson-Lee-Smajlovic (\cite{jorgenson2024functional}),
Malo-Zung (\cite{2024zetafunctionsfriedconjecture}) and
Humbert-Tao (\cite{humbert2026twistedpollicottruelleresonanceszeta}).
We mention also here a related work by 
Abdurrahman-Venkatesh (\cite{abdurrahman2025symplectic}).

\begin{ack}

The author would like to warmly thank the organizers of the thematic trimester programme:
"Representation Theory and Noncommutative Geometry"
for the opportunity to participate and deliver the mini-course:
"Twisted Ruelle zeta function, complex-valued analytic torsion and the Fried's conjecture".
The author was supported by the Mottez Fellowship.

\end{ack}

\section{Dynamical zeta functions}

In Subsections \ref{lss1} and \ref{tdzf}, we recall some basic facts 
and definitions for locally symmetric spaces
of real rank one and the twisted dynamical zeta functions,
mainly from Bunke-Olbrich \cite{BO}.

\subsection{Locally symmetric spaces of real rank one}
\label{lss1}

Let $X$ be a compact locally symmetric Riemannian manifold of negative sectional curvature. 
The universal covering $\widetilde{X}$ of $X$ is a Riemannian symmetric space of real rank one.
Then, $\widetilde{X}$ is a real, complex or quaternionic hyperbolic space or the hyperbolic Cayley space. More precisely, 
\begin{align*}
&\widetilde{X}= G/ K,\\
&X=\Gamma\backslash G/ K
\end{align*}
where $G$ is the quotient of one of the groups 
\begin{equation*}
\Spin(d,1), \SU(d,1), \Sp(d), \F_{4}^{-20},
\end{equation*}
by a subgroup of its center,
$K$ is a maximal subgroup of $G$ and $\Gamma$
is a discrete, co-compact, torsion-free subgroup of $G$.

Let $\mathfrak{g}=\mathfrak{g}\oplus\mathfrak{p}$
be the Cartan decomposition of the Lie algebra $\mathfrak{g}$ of $G$.
Let $\mathfrak{a}$ be a maximal abelian subalgebra of $\mathfrak{p}$.
Let $\Phi(\mathfrak{g},\mathfrak{a})$ be the root system determined by
the adjoint action of $\mathfrak{a}$ on $\mathfrak{g}$.
Let $\Phi^{+}(\mathfrak{g},\mathfrak{a})\subset \Phi(\mathfrak{g},\mathfrak{a})$
be a fixed system of positive roots and $\mathfrak{n}_{\alpha}$ the root spaces,
corresponding to $\alpha\in \Phi^{+}(\mathfrak{g},\mathfrak{a})$.
Let
\begin{equation*}
\mathfrak{n}=\sum_{\alpha\in \Phi^{+}(\mathfrak{g},\mathfrak{a})}\mathfrak{n}_{\alpha}.
\end{equation*}
Then, the Iwasawa decomposition of $\mathfrak{g}$ is 
\begin{equation*}
\mathfrak{g}=\mathfrak{k}\oplus\mathfrak{a}\oplus\mathfrak{n},
\end{equation*}
which corresponds to the Iwasawa decomposition of $G$
\begin{equation*}
G=KAN.
\end{equation*}
$X$ has real rank one, meaning that 
\begin{equation*}
\dim(\mathfrak{a})=1.
\end{equation*}

We define $\rho\in\mathfrak{a}_{\C}^{*}$ by
\begin{equation}\label{rho}
\rho:=\frac{1}{2}\sum_{\alpha\in \Phi^{+}(\mathfrak{g},\mathfrak{a})} 
\dim\mathfrak{n}_{\alpha }\alpha.
\end{equation}

The positive Weyl chamber $\mathfrak{a}^{+}$ is defined as the half line 
in $\mathfrak{a}$, on which the positive roots take 
positive values. Let $A^{+}:=\exp(\mathfrak{a}^{+})\subset A$.

From geometry point of view, the rank of $\widetilde{X}$ is 
the dimension of the maximal flat subspace of 
$\widetilde{X}$.
In our case, that is, when the rank of $\widetilde{X}$ is one, 
the \textit{maximal flat subspaces} are the \textit{geodesics}.

Let $e$ be the identity element of $G$.
Recall that there is a canonical isomorphism 
\begin{equation*}
 T_{eK}\widetilde{X}\cong \mathfrak{g}/\mathfrak{k}
 \cong\mathfrak{p} .
\end{equation*}
We get a metric on $\mathfrak{p}$,
by restricting to $\mathfrak{p}\times\mathfrak{p}$
the normalized $\Ad(G)$-invariant inner product 
on $\mathfrak{g}\times\mathfrak{g}$.
Let $M$ be the centralizer of $\mathfrak{a}$ in $K$ 
with Lie algebra $\mathfrak{m}$.
Since the subgroup $K$ acts transitively on the
unit vectors in $\mathfrak{p}$, 
we have for the unit sphere bundle of $\widetilde{X}$,
\begin{equation*}
S(\widetilde{X})=G/M,
\end{equation*}
and correspondingly for the unit sphere bundle of $X$ 
\begin{equation*}
S(X)=\Gamma\backslash G/M.
\end{equation*}

\subsection{Twisted dynamical zeta functions-Definitions}
\label{tdzf}

We want to define the twisted Selberg and Ruelle zeta functions,
associated with the \textit{geodesic flow} on $S(X)$.
We know that in hyperbolic geometry there is one-to-one 
correspondence between the set of closed geodesics $c_{\gamma}$ on the hyperbolic manifold 
$X$ and the set of conjugacy classes $[\gamma]$ in $\pi_{1}(X)$.
Among the closed geodesics on $X$, we consider the 
\textit{prime} ones, i.e., those closed geodesics
that they trace their image exactly once and hence
their \textit{winding number} is \textit{one}.
Recall from \cite[Lemma 6.5]{Wa}
that in our setting \ref{lss1}, every element in 
$\Gamma$, apart from the identity, is conjugate in $G$ to an element in $MA^{+}$.

We recall here the definitions of the \textit{twisted Selberg and Ruelle zeta functions}
from Bunke and Olbrich (\cite{BO}).
Let $\sigma$ and $\chi$ be finite-dimensional, \textit{unitary representations} of 
of $M$ and $\Gamma$, correspondingly. 

\begin{defi}(\cite[Definition 3.2]{BO})

For $s\in \C$ with $\RE(s)>\rho$, the twisted Selberg zeta function is defined by
\begin{align*}
Z(s;\sigma,\chi):=\prod_{\substack{[\gamma]\neq{e}\\ \\ [\gamma]\prim}} &\prod_{k=0}^{\infty}\det\bigg
(\Id-(\chi(\gamma)\otimes\sigma(m_\gamma)\otimes \\
&S^k(\Ad(m_\gamma a_\gamma)|_{\overline{\mathfrak{n}}}))
e^{-(s+|\rho|)\lvert {l(\gamma)}}\bigg),
\end{align*}
where $\overline{\mathfrak{n}}$ is the sum of the negative root spaces of $\mathfrak{a}$,
$S^k(\Ad(m_\gamma a_\gamma)_{\overline{\mathfrak{n}}})$ denotes the $k$-th
symmetric power of the adjoint map $\Ad(m_\gamma a_\gamma)$ restricted to $\mathfrak{\overline{n}}$ ($\rho$ is as in \eqref{rho}).
\end{defi}

\begin{defi}(\cite[Definition 3.1]{BO})

For $s\in \C$ with $\RE(s)>2\rho$, the twisted Ruelle zeta function is defined by
\begin{equation*}
 R(s;\sigma,\chi):=\prod_{\substack{[\gamma]\neq{e}\\ \\ [\gamma]\prim}}\det\big(\Id-(\chi(\gamma)
 \otimes\sigma(m_{\gamma}))e^{-s{{l(\gamma)}}}\big)^{(-1)^{d-1}}
\end{equation*}
($\rho$ is as in \eqref{rho}).
\end{defi}

\subsection{Non-unitary representations}

From now on, we consider arbitrary, complex, finite-dimensional representation of $\Gamma$
\begin{equation*}
\chi\colon\Gamma\rightarrow \GL(V_{\chi}).
\end{equation*}

The first problem that one has to encounter is the domain of convergence for the 
of the twisted dynamical zeta functions. 
Since the representation is no longer unitary, one has to 
study possible bounds for the trace $  \lvert\tr(\chi(\gamma))\rvert$.

\subsubsection{Real hyperbolic manifolds of odd dimension}
We consider compact real hyperbolic manifolds $X$ of odd dimension $d$.
By \cite[Lemma 3.3]{Spilioti2018}, we have the following lemma. 
\begin{lem}\label{lemmarep}
 Let $\chi\colon\Gamma\rightarrow \GL(V_{\chi})$ be a finite-dimensional representation of $\Gamma$.
 Then, there exist positive constants $\alpha,\beta>0$ such that 
 \begin{equation}\label{cruciallemma}
  \lvert\tr(\chi(\gamma))\rvert\leq \alpha e^{\beta l(\gamma)},\quad \forall\gamma\in\Gamma-\{e\}.
 \end{equation}
\end{lem}

Then, we use the fact that for $\Gamma$ a cocompact lattice, 
there exists a positive constant $C'$ such that 
\begin{equation}\label{geodest}
\sharp \{[\gamma]:l(\gamma)<R\}\leq\sharp\{ \gamma \in \Gamma: l(\gamma)\leq R\}\leq C'e^{2|\rho |R}
\end{equation}
(\cite[(1.31)]{BO}).

Using \eqref{cruciallemma} and \eqref{geodest},
we can define the \textit{twisted dynamical zeta functions for non-unitary representations} $\chi$
(for more details see \cite[Section 3]{Spilioti2018}).

\begin{defi}
The twisted Selberg zeta function $ Z(s;\sigma,\chi)$ is defined for $Re(s)>c_{1}$ by
 \begin{equation}
 Z(s;\sigma,\chi):=\prod_{\substack{[\gamma]\neq{e}\\ [\gamma]\prim}}\prod_{k=0}^{\infty}\det(\Id-(\chi(\gamma)\otimes\sigma(m_\gamma)\otimes S^k(\Ad(m_\gamma a_\gamma)_{\overline{\mathfrak{n}}}))e^{-(s+|\rho|)l(\gamma)}),
\end{equation}
where $c_{1}$ is a positive constant.
\end{defi}

\begin{defi}\label{ruelledefiodd}
The twisted Ruelle zeta function $ R(s;\sigma,\chi)$ is defined for $\Re(s)>c_{2}$ by
 \begin{equation}
 R(s;\sigma,\chi):=\prod_{\substack{[\gamma]\neq{e}\\ [\gamma]\prim}}\det\big(\Id-(\chi(\gamma)\otimes\sigma(m_{\gamma}))e^{-sl(\gamma)}\big),
\end{equation}
where $c_{2}$ is a positive constant.
\end{defi}

\subsubsection{Compact hyperbolic surfaces}
Let now $X$ be  compact hyperbolic surface.
To define the twisted dynamical zeta functions for \textit{non-unitary representations} $\chi$,
we use the fact that the word length of an element in $\Gamma$ is quasi isometric to the hyperbolic distance (see \cite[Proposition 3.2]{Lmr} and \cite[p.10]{Wo}).
We have that there exists a $c\geq 0$ such that 
\begin{equation}\label{chibound}
	\Vert\chi(\gamma)\Vert\leq e^{cl(\gamma)},
	\end{equation}
where $\Vert\cdot\Vert$ is the operator norm associated to a fixed norm on $V_{\chi}$.
We use also estimates \eqref{geodest}.
Then, we have the following definitions
(for more details see \cite[Section 4.1]{frahm2023twisted}).

\begin{defi}
The twisted Selberg zeta function is defined for $\Re(s)>c_{3}$ by
\begin{equation}
		Z(s;\chi):=\prod_{\substack{[\gamma]\neq{e}\\ [\gamma]\prim}} \prod_{k=0}^{\infty}\det\big(\Id-\chi(\gamma)e^{-(s+k)l(\gamma)}\big),
	\end{equation}
where $c_{3}$ is a positive constant.
\end{defi}

\begin{defi}
The twisted Ruelle zeta function is defined for $\Re(s)>c_{4}$ by
\begin{equation}
		R(s;\chi)=\prod_{\substack{[\gamma]\neq{e}\\ [\gamma]\prim}}\det (\Id-\chi(\gamma)e^{-sl(\gamma)}),
	\end{equation}
where $c_{4}$ is a positive constant.
\end{defi}

\begin{rmrk}
\begin{enumerate}
\item The constants $c_{3}, c_{4}$ can be explicitly determined, in terms of 
the constant $c$ as in \eqref{chibound}. 
One can consider the notion of the \textit{critical exponent} of $
\chi$, in order to eliminate the dependence on the norm $\Vert \cdot\Vert$ on $V_{\chi}$
in the estimates \eqref{chibound}. The critical exponent $c_{\chi}$ is defined as
$c_{\chi}:=\inf\{c\geq 0, \text{such that \eqref{chibound} holds}\}$.

\item In \cite{naud2023determinants}, Naud and the author, used the notion 
of the critical exponent to study the \textit{bottom of the spectrum of the twisted 
Laplacian}, and for some families of linear representations describe it more 
explicitly.

\item The constants $c_{1}$ and $c_{2}$ can be also determined explicitly in terms of 
the term $|\rho|$, and the constant $\beta$ in \eqref{cruciallemma}.
In fact, Lemma \ref{lemmarep} is also a consequence of the result of 
Lubotzky, Mozes and Raghunathan (\cite[Proposition 3.2]{Lmr}).

\end{enumerate}
\end{rmrk}

\subsubsection{Compact hyperbolic orbisurfaces}

We briefly recall some basic facts about hyperbolic obirsurfaces.
For more details,  we refer the reader to \cite[Section 2]{fried1986fuchsian} and
\cite[Section 1]{JEP_2023__10__1391_0}.
Let $G=\PSL(2,\R)$ and $K=\PSO(2)$.
Let $\Gamma\subseteq G$ be a cocompact Fuchsian group, meaning a discrete subgroup of $G$ such that the quotient $X=\Gamma\backslash \H^2=\Gamma\backslash G/K$ is compact.  The Poincar\'{e} metric on~$\H^2$ induces a metric on $X$ that turns it into a compact hyperbolic \textit{orbisurface}. 
In this case, every \textit{non-trivial} element in $\Gamma$ is either \textit{hyperbolic or elliptic}.

Let $\widetilde{G}$ denote the universal cover of $G$ and write~$\widetilde{H}$ for the preimage of a subgroup $H\subseteq G$ under the universal covering map. Consider the following one-parameter family in $G$ (modulo $\pm I_2$):
$$
k_\theta=\begin{pmatrix}\cos\theta&\sin\theta\\-\sin\theta&\cos\theta\end{pmatrix},
$$
where $\theta\in\R$.
Denote by $\widetilde{k}_\theta$ the unique lift to the universal cover $\widetilde{G}$, turning
$$ \widetilde{K}=\{\widetilde{k}_\theta:\theta\in\R\}$$
into a connected one-parameter subgroup of $\widetilde{G}$.
The center $\widetilde{Z}$ of $\widetilde{G}$ is given by
$ \widetilde{Z} = \{\widetilde{k}_\theta:\theta\in\Z\pi\}. $

The unit tangent bundle of $\H^2$ naturally identifies with $G=\PSL(2,\R)$. Indeed, the map
$$
G \to \H^2
\,\,\,\,\,\textrm{\rm given by} \,\,\,\,\,
g \mapsto g\cdot \textrm i
$$
is a fibration whose fiber is a circle $S^1$ of directions of unit tangent vectors.
The quotient $X_{1}:=\Gamma \backslash G$ is a \textit{compact 3-manifold}, which we will call the unit tangent bundle of $X$.
It follows that $X_{1}$ is a \textit{Seifert fibered 3-dimensional manifold}, and this fibration induces the exact sequence
$$
1 \to \Z = \pi_1(\PSO(2)) \to \pi_1(X_1) \to \Gamma= \pi_1(X) \to 1
$$
where $\pi_1(X)$ denotes the orbifold fundamental group of $X$.

To define the twisted dynamical zeta functions, we note first that 
every closed oriented geodesic $\gamma$ on $X$ lifts canonically to a homotopy class in $\widetilde{\Gamma}=\pi_1(X_1)$,
which we also denote by $\gamma$.  
The conjugacy class $[\gamma]$ in $\widetilde{\Gamma}/\widetilde{Z}\simeq\Gamma=\pi_1(X)$
consists of hyperbolic elements.  Conversely, every such conjugacy class of hyperbolic elements contains exactly one representative of a closed geodesic.

\begin{defi}
Let $\rho:\widetilde{\Gamma}=\pi_1(X_1)\to\GL(V_\rho)$ be a finite-dimensional, complex representation. 
For $s\in\C$, we define the twisted Selberg zeta function by
\begin{equation*}
Z(s;\rho) = \prod_{\gamma\text{ prime}}\prod_{k=0}^\infty\det\left(\Id-\rho(\gamma)e^{-(s+k)\ell(\gamma)}\right).
\end{equation*}
For $s\in\C$, we define the twisted Ruelle zeta function by
\begin{equation*}
R(s;\rho) = \prod_{\gamma\text{ prime}}\det\left(\Id-\rho(\gamma)e^{-s\ell(\gamma)}\right).
\end{equation*}
\end{defi}
It is shown in \cite[Theorem 3.1]{fedosova2020meromorphic} that $Z(s;\rho)$, and hence $R(s;\rho)$, converges for $s$ in some right half plane of $\C$ and defines a holomorphic function on this half plane (see also \cite[Section 1.2]{Wo} for the torsion-free case).

\section{Trace formulas}\label{tracecases}

The \textit{Selberg trace formula} is a powerful tool from 
\textit{harmonic analysis on symmetric spaces} that can be used 
to study the dynamical zeta functions. Under different geometrical 
and algebraic settings, it can be a complicated to derive and use
the trace formula. The advantage is that one can study all the locally symmetric spaces, 
using the same method.

The \textit{Selberg trace formula for non-unitary representations} was developed by 
M{\"u}ller in \cite{M}. It is derived for all locally symmetric spaces $X$ of real rank one and
Paley-Wiener test functions. In this paper, M{\"u}ller  considered the
\textit{twisted Laplacians}, which are elliptic, \textit{non-self-adjoint operators}.

We recall here briefly the definition.
Let  $E_{\chi}\rightarrow X$ be the flat vector bundle associated with $\chi$
and  $C^{\infty}(X,E_{\chi})$ be 
the space of smooth sections of the
flat vector bundle $E_{\chi}$.
Let $\ast_{\chi}$ be the extension of the Hodge star operator $\ast$ to 
an operator in $\Lambda^{p}(T^{*}X)\otimes E_{\chi}$, defined by
\begin{equation*}
\ast_{\chi}=\ast\otimes\Id_{E_{\chi}}.
\end{equation*}
We define
\begin{equation*}
\delta^{\sharp}_{\chi}:=(-1)^{n+1}\ast_{\chi} d_{\chi}\ast_{\chi},
\end{equation*}
and 
\begin{equation*}
\Delta^{\sharp}_{\chi}=\delta^{\sharp}_{\chi}d_{\chi},
\end{equation*}
where $d_{\chi}$ is the associated exterior derivative.

One can chose a Hermitian metric on $E_{\chi}$ and use it 
together with the Riemannian metric on $X$ to introduce an inner product on
$C^{\infty}(X,E_{\chi})$. Nevertheless, the twisted Laplacian is \textit{not} in general a 
self-adjoint operator with respect to this inner product.
Let $\delta_{\chi}$  be the formal adjoint of $d_{\chi}$ with respect 
to the inner product as above.
The operator $\delta^{\sharp}_{\chi}$ is \textit{not} the formal adjoint of 
$d_{\chi}$ with respect to the inner product and is rather a perturbation of 
$\delta_{chi}$ by a zero order operator (see \cite[p. 177-178]{cm}).
Let $\Delta_{\chi}$ be the associated self-adjoint operator 
$\Delta_{\chi}=\delta_{\chi}d_{\chi}$.
The twisted Laplacian $\Delta_{\chi}^{\sharp}$ has the same principal symbol as 
$\Delta_{\chi}$ and hence has \textit{nice spectral properties}. 
Its spectrum is discrete and contained in a translate of a positive cone 
$C \subset\C$ with $R^{+}\subset C$ (\cite{Sh}).

A Selberg trace formula for non-unitary twists and the 
heat operator $e^{-t\Delta^{\sharp}_{\chi}}$ is derived 
in \cite{Spilioti2018} for compact hyperbolic odd-dimensional manifolds,
in \cite{frahm2023twisted} for compact hyperbolic surfaces and in 
in \cite{frahm2023twisted} for compact hyperbolic orbisurfaces.

The  trace formulas \eqref{selberg1}, \eqref{selberg2}, 
\eqref{selberg3a} and \eqref{selberg3b}
consists of two sides: the 
\textit{spectral} (right-hand side) and the \textit{geometrical one} (left-hand side).
We prove that the heat operator induced by the twisted Laplacian is 
an integral operator with smooth kernel and of trace class.
Then we consider the trace of this operator, which 
can be explicitly written by Lidskii's theorem and 
also develop the integrals that involve the trace of the \textit{heat kernel}.

We provide here some more details in the case of the \textit{compact
hyperbolic surfaces}, but the general scheme holds for
the case of a compact hyperbolic odd-dimensional manifolds
and a hyperbolic orbisurface as well.
In the last two cases the representation theory
of the the Lie groups involved can be more complicated, 
for more details see \cite{Spilioti2018} and \cite{JEP_2023__10__1391_0}.

For the derivation of the trace formula, we follow the work of Wallach (\cite{Wab}), and its extension to non-unitary representations setting by M\"{u}ller (\cite{M}).
We denote by $\spec(\Delta^{\sharp}_{\chi})\subseteq\C$ the discrete) spectrum of $\Delta^{\sharp}_{\chi}$. For $\mu\in \spec(\Delta^{\sharp}_{\chi})$, we write $L^2(X,E_\chi)_\mu$ for the corresponding generalized eigenspace. We define the algebraic multiplicity $m(\mu)$ of $\mu$ as $m(\mu):=\dim L^2(X,E_\chi)_\mu$. We use the extended results for the heat operator
obtained in \cite[p. 171--173]{Spilioti2018} to conclude that the heat operator $e^{-t{\Delta}^{\sharp}_{\chi}}$ is an integral operator with smooth kernel, i.e., there exists a smooth section $H^{\chi}_{t}$ of $\End(E_\chi)$ on $X\times X$ such that for $f\in L^{2}(X,E_{\chi})$, we have
	\begin{equation*}
		e^{-t{\Delta}^{\sharp}_{\chi}}f(x)=\int_{X}H^{\chi}_{t}(x,y)f(y)dy.
	\end{equation*}
	By \cite[Proposition 2.5]{M}, this operator is of trace class.
	By Lidskii's theorem~\cite[Theorem 3.7]{SB}, we have 
	\begin{equation*}
		\tr(e^{-t{\Delta}^{\sharp}_{\chi}})=\sum_{\mu\in
			\spec(\Delta^{\sharp}_{\chi})}m(\mu)e^{-t\mu},
	\end{equation*}
	which is the spectral side of our (pre)-trace formula.
	Let $H_{t}$ be the kernel of the heat operator $e^{-t\widetilde{\Delta}}$, that is the operator induced by the self-adjoint Laplacian $\widetilde{\Delta}$ acting in $L^2(\widetilde{X})$. By \cite[Lemma 2.3 and Proposition 2.4]{BM}, $H_t$ is contained in the Harish-Chandra $L^q$-Schwartz space $\mathcal{C}^q(G)$ for any $q>0$ (see e.g. \cite[p. 161--162]{BM} for the definition of the Schwartz space). 
Moreover, by \cite[eq. (3.7)]{frahm2023twisted} we have
	\begin{equation*}
		H_{t}^\chi(x,x')=\sum_{\gamma \in \Gamma}H_{t}(g^{-1}\gamma g')\chi(\gamma),
	\end{equation*}
	where $x=\Gamma g K,x'=\Gamma g' K\in X$ with $g,g'\in G$.
 By \cite[Proposition 5.3]{Spilioti2018}, we have the following result.
	
	\begin{prop}
		Let $E_{\chi}$ be a flat vector bundle over $X=\Gamma\backslash \widetilde{X}$ associated with a finite-dimensional, complex
		representation $\chi\colon\Gamma\rightarrow \GL(V_{\chi})$ of $\Gamma$. Let $\Delta_{\chi}^{\sharp}$ be the twisted Bochner--Laplace operator acting in $L^2(X,E_{\chi})$. Then, 
		\begin{equation*}
			\sum_{\mu\in\spec(\Delta_{\chi}^\sharp)}m(\mu)e^{-t\mu} = \tr(e^{-t\Delta_{\chi}^{\sharp}}) =\int_{\Gamma\backslash G}\left(\sum_{\gamma \in \Gamma}\tr \chi(\gamma)\cdot H_{t}(g^{-1}\gamma g)\right)\,d\dot{g}.
		\end{equation*}
	\end{prop}
	
	As in \cite[Proposition 5.1, Proposition 6.1]{M} and \cite[p. 172--173, 177--178]{Wa}, we group the summation into conjugacy classes and use the Fourier inversion formula to obtain
	\begin{multline}\label{tracewallach}
		\tr(e^{-t\Delta_{\chi}^{\sharp}})=\dim(V_{\chi})\Vol(X)H_{t}(e)\\
		+\sum_{[\gamma]\neq e} \tr\chi(\gamma)\frac{l(\gamma)}{n_{\Gamma}(\gamma)D(\gamma)}\frac{1}{2\pi}
		\int_{\R}\Theta_\lambda(H_{t})e^{-il(\gamma)\lambda}d\lambda,
	\end{multline}
	where
	$$ D(\gamma) = e^{-\frac{l(\gamma)}{2}} |\det(\Ad(a_{\gamma})|_{\mathfrak{n}}-\Id)|. $$
Every $\gamma\in\Gamma$, $\gamma\neq e$, is hyperbolic, so the summation is over hyperbolic conjugacy classes of $\Gamma$. We refer to the first term in the right-hand side of 
	 of \eqref{tracewallach}
	as the \emph{identity contribution} and to the second term as the \emph{hyperbolic contribution}.

The \textit{character} $\Theta_\lambda$ of $\pi_\lambda$ can be evaluated on the $K$-biinvariant function $H_t$ in terms of the \textit{spherical Fourier transform} $\widetilde{H_{t}}(\lambda)$ of $H_{t}$:
	\begin{equation*}
		\Theta_{\lambda}(H_{t}) = \tr\pi_{\lambda}(H_{t}) = \int_{G} \langle\pi_{\lambda}(g)v,v\rangle H_{t}(g)dg = \int_{G} \phi_\lambda(g)H_{t}(g)dg = \widetilde{H_{t}}(\lambda),
	\end{equation*}
	where $v\in \mathcal{H}_\lambda^K$, with $\lVert v \rVert=1$ and $\phi_\lambda(g)=\langle\pi_{\lambda}(g)v,v\rangle$ denotes the associated spherical function (see e.g. \cite[Chapter IV]{Hel} for details).

Be \cite[Lemma 3.2.3]{frahm2023twisted}, we have 
\begin{equation*}
			\Theta_{\lambda}(H_{t})=e^{-t(\lambda^2+\frac{1}{4})}.
		\end{equation*}
	
	For the identity contribution in \eqref{tracewallach}, we have that since the
function $H_{t}$ on $G$ belongs to the Harish-Chandra $L^{q}$-Schwartz space, the Fourier inversion formula \cite[Theorem 3]{HC} can be applied to $H_{t}$ (see also Theorem 1 in  \cite{ANKER1991331} for $K$-biinvariant test functions in the Harish-Chandra $L^{q}$-Schwartz space for $0<q\leq2$).
	By \cite[Theorem~7.5~(i)]{Hel} (see also \cite[eq. (28) and (29), p. 42]{Hel}), we have
	\begin{equation*}
		H_{t}(e)=\frac{1}{4\pi^2}\int_{\R}\Theta_{\lambda}(H_{t})\rvert c(\lambda)\lvert^{-2} d\lambda
	\end{equation*}
	with the $c$-function $c(\lambda)$ given by (note the change of variables $\lambda$ vs. $\frac{\lambda}{2}$ compared to \cite{Hel})
	\begin{equation*}
		\lvert c(\lambda)\lvert^{-2}=\lambda\pi\tanh\lambda\pi,\quad \lambda\in\R.
	\end{equation*}
	
From the above discussion we have the following theorem.		
	\begin{thm}(Trace formula, \cite[Theorem 3.2.4]{frahm2023twisted})
		Let $\chi\colon\Gamma\rightarrow\GL(V_{\chi})$ be a finite-dimensional representation of $\Gamma$. 
		Then, the following Selberg trace formula for the operator $e^{-t\Delta_{\chi}^{\sharp}}$ holds:
		\begin{multline}\label{selberg1}
			\tr(e^{-t\Delta_{\chi}^{\sharp}})=\frac{1}{4\pi^2}\dim(V_{\chi})\Vol(X)\int_{\R}e^{-t(\lambda^2+\frac{1}{4})}\lambda\pi\tanh\lambda\pi d\lambda\\
			+\sum_{[\gamma]\neq e} \tr\chi(\gamma)\frac{l(\gamma)}{n_{\Gamma}(\gamma)D(\gamma)}\frac{1}{2\pi}
			\int_{\R}e^{-t(\lambda^2+\frac{1}{4})}e^{-il(\gamma)\lambda}d\lambda.
		\end{multline}
	\end{thm}

\subsubsection{Compact hyperbolic orbisurfaces}

For a compact hyperbolic \textit{orbisurface}, the Selberg trace formula with non-unitary twists is given by Theorem \ref{traceorbi}.
We consider the unitary characters $\tau_m$, $m\in\R$, defined by
\begin{equation}
	\tau=\tau_m(\widetilde{k}_\theta) = e^{im\theta} \qquad (\theta\in\R).
\end{equation}

Let $\rho:\pi_1(X_1)\to\GL(V_\rho)$ be a finite-dimensional, \textit{arbtirary}, complex representation of $\pi_1(X_1)$. For  $u=\widetilde{k}_\pi$, $\rho(u)=e^{-i\pi m} {\rm Id}_{V_\rho}$ for some rational number $m$ (see \cite{JEP_2023__10__1391_0}, p. 1417).
The operator $A_{\tau_m,\rho}^{\sharp}$ is a 
\textit{shift} of the twisted Laplacian
$\Delta_{\tau_m,\rho}^{\sharp}$ (see \cite[Section 3.7]{JEP_2023__10__1391_0}).

\begin{thm}(Trace formula, \cite[Theorem 3.7.1]{JEP_2023__10__1391_0})
\label{traceorbi}
For every $t>0$ we have
	\begin{align}\label{selberg2}
		& \Tr(e^{-tA_{\tau_m,\rho}^\sharp})\\\notag
		& \quad=\frac{\Vol(X)\dim(V_\rho)}{4\pi}\Bigg[\int_\R e^{-t\lambda^2}\frac{\lambda\sinh(2\pi\lambda)}{\cosh(2\pi\lambda)+\cos(\pi m)}\,d\lambda\\\notag
		& \hspace{7cm}+ \sum_{\substack{1\leq\ell<|m|\\\ell\textup{ odd}}}(|m|-\ell)e^{(\frac{|m|-\ell}{2})^2t}\Bigg].\\\notag
		& \qquad+\frac{1}{2\sqrt{4\pi t}}\sum_{[\gamma]\textup{ hyp.}}\frac{\ell(\gamma)\tr\rho(\gamma)}{n_\Gamma(\gamma)\sinh(\frac{\ell(\gamma)}{2})}e^{-\frac{\ell(\gamma)^2}{4t}}\\\notag
		& \qquad+\sum_{[\gamma]\textup{ ell.}} \frac{\tr\rho(\gamma)}{4M(\gamma)\sin(\theta(\gamma))}\Bigg[\int_\R e^{-t\lambda^2}\frac{\cosh(2(\pi-\theta(\gamma))\lambda)+e^{i\pi m}\cosh(2\theta(\gamma)\lambda)}{\cosh(2\pi\lambda)+\cos(\pi m)}\,d\lambda\\
		& \hspace{4cm}+2i\sign(m)\sum_{\substack{1\leq\ell<|m|\\\ell\textup{ odd}}}e^{i\sign(m)(|m|-\ell)\theta(\gamma)}e^{(\frac{|m|-\ell}{2})^2t}\Bigg].
	\end{align}
	Here, the summation is over the conjugacy classes $[\gamma]$ of hyperbolic, resp. elliptic, elements in $\widetilde{\Gamma}/\widetilde{Z}\simeq\Gamma$.
\end{thm}

\subsubsection{Compact hyperbolic odd-dimensional manifolds}

The Selberg trace formula with non-unitary twists in the case of a compact hyperbolic 
\textit{odd-dimensional} manifold $X=\Gamma\backslash \H^{d}$ is derived in 
\cite{Spilioti2018}.

Let $\chi\colon\Gamma\rightarrow \GL(V_{\chi})$ be
a finite-dimensional representation of $\Gamma$.
Let  $M'=\Norm_{K}(A)$ be the normalizer of $A$ in $K$.
We define the restricted Weyl group as the quotient $ W_{A}:=M'/M$
(see Section \ref{lss1} for the definition of the groups $A$, $K$ and $M$).
Let $w\in W_{A}$ be the non-trivial element of $W_{A}$, and 
$m_{w}$ a representative of $w$ in $M'$.
The action of $W_{A}$ on $\widehat{M}$ is defined by
\begin{equation*}
(w\sigma)(m):=\sigma(m_{w}^{-1}mm_{w}),\quad m\in M, \sigma\in\widehat{M}.
\end{equation*}

We have the following theorem.

\begin{thm}(Trace formula, \cite[Theorem 5.5]{Spilioti2018})
 Let $\sigma \in \widehat{M}$.
 \begin{itemize}
 \item If $\sigma$ is Weyl invariant, 
 \begin{align} \label{selberg3a}
 \begin{aligned}
  \Tr(e^{-tA_{\chi}^{\sharp}(\sigma)})=&\dim(V_{\chi})\Vol(X)\int_{\R}e^{-t\lambda^{2}}P_{\sigma}(i\lambda)d\lambda\\
  &+\sum_{[\gamma]\neq e} \frac{l(\gamma)}{n_{\Gamma}(\gamma)}L(\gamma;\sigma,\chi)
  \frac{e^{-l(\gamma)^{{2}}/4t}}{(4\pi t)^{1/2}};
  \end{aligned}
 \end{align}
  \item If $\sigma$ is non-Weyl invariant, 
   \begin{align}\label{selberg3b}
   \begin{aligned}
  \Tr(e^{-tA_{\chi}^{\sharp}(\sigma)})=&2\dim(V_{\chi})\Vol(X)\int_{\R}e^{-t\lambda^{2}}P_{\sigma}(i\lambda)d\lambda\\
  &+\sum_{[\gamma]\neq e}\frac{l(\gamma)}{n_{\Gamma}(\gamma)}L(\gamma;\sigma+w\sigma,\chi)
  \frac{e^{-l(\gamma)^{{2}}/4t}}{(4\pi t)^{1/2}},
  \end{aligned}
 \end{align}
\end{itemize} 
where  
 \begin{equation*}
 L(\gamma;\sigma,\chi)= \frac{\tr(\chi(\gamma)\otimes\sigma(m_{\gamma}))e^{-|\rho|l(\gamma)}}{\det(\Id-\Ad(m_{\gamma}a_{\gamma})_{\overline{n}})}.
\end{equation*}
\end{thm}

Here, $A_{\chi}^{\sharp}(\sigma)$ is a certain elliptic operator, induced by the twisted,
non-self-adjoint Laplacian $\Delta^{\sharp}_{\chi}$.

\section{Special value at zero}

Having the trace formula at hand, one can prove certain properties of the dynamical zeta functions
in the cases as in Section \ref{tracecases}.
The key-idea is that in all these cases, after certain technicalities, the \textit{hyperbolic contribution in the geometrical part of the trace formula} is related to the \textit{logarithmic derivative of the 
twisted Selberg zeta function}. 

We provide here some more details in the case of a compact hyperbolic surface.
For the cases of compact hyperbolic orbisurfaces and odd-dimensional  compact 
hyperbolic manifolds, we refer the reader to \cite{JEP_2023__10__1391_0} and \cite{Spilioti2018}.

\subsection{Idea of the proof}
\subsubsection{Meromorphic continuation and functional equations}

	The trace formula \eqref{selberg1} motivates to consider a shift of the operator $\Delta_{\chi}^{\sharp}$ by $\frac{1}{4}$:
	\begin{equation*}
		A^{\sharp}_{\chi}:=\Delta_{\chi}^{\sharp}-\frac{1}{4}.
	\end{equation*}

By \cite[Corollary 4.2.1]{frahm2023twisted}, we have the following corollary.
	\begin{coro}
		Let $X=\Gamma\backslash\H^2$ be a compact hyperbolic surface and $\chi$ be a finite-dimensional, complex representation of $\Gamma$. Then, we have
		\begin{multline}\label{selberg2a}
			\tr(e^{-tA^{\sharp}_{\chi}})=\frac{1}{4\pi^2}\dim(V_{\chi})\Vol(X)\int_{\R}e^{-t\lambda^2}\lambda\pi \tanh\lambda\pi d\lambda\\
			+\frac{1}{2\sqrt{4\pi t}}\sum_{[\gamma]\neq e} \frac{l(\gamma)\tr(\chi(\gamma))}{n_{\Gamma}(\gamma)\sinh(l(\gamma)/2)}
			e^{-\frac{l(\gamma)^{2}}{4t}}.
		\end{multline}
	\end{coro}

We want to use a \textit{resolvent trace formula}.	
Let $s_{1},s_{2}\in\C$ such that $(s_{1}-1/2)^2,(s_{2}-1/2)^2\in \C  \setminus\spec(-A^{\sharp}_{\chi})$.
We consider the product of the two resolvents $(A_\chi^\sharp+(s_j-\frac{1}{2})^2)^{-1}$, for $j=1,2$.
By \cite[Lemma 2.2]{M}, the operator $A_\chi^\sharp$ satisfies \textit{Weyl's Law}$A_\chi^\sharp$. This implies that
	\begin{equation*}
		\big(A^{\sharp}_{\chi}+(s_{1}-1/2)^2\big)^{-1}\big(A^{\sharp}_{\chi}+(s_{2}-1/2)^2\big)^{-1}
	\end{equation*}
	is of trace class, and so is
	\begin{multline}\label{resid}
		\big(A^{\sharp}_{\chi}+(s_{1}-1/2)^2\big)^{-1}-(A^{\sharp}_{\chi}+(s_{2}-1/2)^2\big)^{-1}\notag\\
		=\big((s_{2}-1/2)^2-(s_{1}-1/2)^2\big)\big(A^{\sharp}_{\chi}+(s_{1}-1/2)^2\big)^{-1}\big(A^{\sharp}_{\chi}+(s_{2}-1/2)^2\big)^{-1}.
	\end{multline}
	We observe also that for $\RE(s_j)\gg0$
	\begin{equation*}
		(A^{\sharp}_{\chi}+(s_j-1/2)^2)^{-1}
		=\int_{0}^{\infty}e^{-t(s_j-1/2)^2}e^{-tA^{\sharp}_{\chi}}\,dt.
	\end{equation*}
	Hence,
	\begin{multline*}
		\big(A^{\sharp}_{\chi}+(s_{1}-1/2)^2\big)^{-1}-\big(A^{\sharp}_{\chi}+(s_{2}-1/2)^2\big)^{-1}\\
		= \int_{0}^{\infty}
		\big(e^{-t(s_{1}-1/2)^2}-e^{-t(s_{2}-1/2)^2}\big)e^{-tA^{\sharp}_{\chi}}dt.
	\end{multline*}
	
Using the \textit{short time asymptotic expansion of the  heat trace} 
(\cite[Lemma 4.2.2]{frahm2023twisted})
\begin{equation}
			\tr e^{-tA^{\sharp}_{\chi}} \sim \sum_{j=0}^{\infty}c_{j}t^{\frac{j-2}{2}} \qquad \mbox{as }
			t\to0^+
		\end{equation}
(for some coefficients $c_{j}$), we obtain the following lemma (\cite[Lemma 4.2.3]{frahm2023twisted}).

	\begin{lem}
		Let $s_{1},s_{2}\in\C$ with $\RE((s_{1}-1/2)^2),\RE((s_{1}-1/2)^2)\gg0$. Then,
		\begin{multline}\label{resolvent}
			\tr\bigg(\big(A^{\sharp}_{\chi}+(s_{1}-1/2)^2\big)^{-1}-\big(A^{\sharp}_{\chi}+(s_{2}-1/2)^2\big)^{-1}\bigg)\\
			= \int_{0}^{\infty}(e^{-t(s_{1}-1/2)^2}-e^{-t(s_{2}-1/2)^2}) \tr(e^{-tA^{\sharp}_{\chi}})dt.
		\end{multline}
	\end{lem}

Let 
\begin{equation*}
L(s;\chi):=\frac{d}{ds}\log Z(s;\chi)
\end{equation*}
be the logarithmic derivative of the twisted Selberg zeta function $Z(s;\chi)$.
By \cite[Lemma 4.1.2]{frahm2023twisted}, it is given by 
		\begin{equation*}
			L(s;\chi)= \sum_{[\gamma]\neq{e}}\frac{l(\gamma)\tr(\chi(\gamma))}{2n_{\Gamma}(\gamma)\sinh(l(\gamma)/2)}e^{-(s-\frac{1}{2})l(\gamma)}.
		\end{equation*}
		
The idea is now to substitute the trace $ \tr(e^{-tA^{\sharp}_{\chi}})$ in the right-hand side of  
\eqref{resolvent} with the \textit{geometrical side of the trace formula} \eqref{selberg2a}.
Then, $L(s;\chi)$ \textit{will appear}.
By \cite[Proposition 4.2.4]{frahm2023twisted}, we have the following proposition.

	\begin{prop}[Resolvent trace formula]\label{resolventtrace}
		Let $s_{1},s_{2}\in\C$ with $\RE(s_1),\RE(s_2)\gg0$ and $\RE((s_{1}-1/2)^2),\RE((s_2-1/2)^2)\gg0$.
		Then, 
		\begin{multline}\label{finalresolvent}
			\tr\bigg(\big(A^{\sharp}_{\chi}+(s_{1}-1/2)^2\big)^{-1}-(A^{\sharp}_{\chi}+(s_{2}-1/2)^2\big)^{-1}\bigg)\\
			=\bigg(\frac{1}{4\pi^2}\dim(V_{\chi})\Vol(X)\bigg)\bigg(\int_{\R}\frac{\lambda\pi\tanh \lambda\pi}{(s_{1}-\frac{1}{2})^{2}+\lambda^{2}}
			-\frac{\lambda\pi\tanh\lambda\pi }{(s_{2}-\frac{1}{2})^{2}+\lambda^{2}}
			d\lambda\bigg)\\
			+\frac{1}{2(s_{1}-\frac{1}{2})}L(s_{1};\chi)
			-\frac{1}{2(s_{2}-\frac{1}{2})}L(s_{2};\chi).
		\end{multline}
	\end{prop}
	
The resolvent trace formula \eqref{finalresolvent} is the \textit{key-formula}
to prove the meromorphic continuation of $L(s;\chi)$.
By \cite[Proposition 4.2.5]{frahm2023twisted}, we have the 
following proposition.

	\begin{prop}\label{merolog}
		The logarithmic derivative $L(s;\chi)$ of the Selberg zeta function $Z(s;\chi)$ extends to a meromorphic function in $s\in\C$ with singularities given by the following formal expression:
		\begin{equation}
			\sum_{j=0}^\infty\left[\frac{1}{s-\frac{1}{2}-i\mu_j}+\frac{1}{s-\frac{1}{2}+i\mu_j}\right] + \frac{\Vol(X)\dim(V_\chi)}{2\pi}\sum_{k=0}^\infty\frac{1+2k}{s+k}.\label{eq:merologFormalSingularities}
		\end{equation}
		where $(\lambda_j=\frac{1}{4}+\mu_j^2)_{j\in\Z^+}\subseteq\C$ are the eigenvalues of $\Delta_\chi^\sharp$ counted with algebraic multiplicity.
	\end{prop}

Using now Proposition~\ref{merolog} and by integration and exponentiation,
we can easily obtain the holomorphic continuation of the twisted Selberg zeta function $Z(s;\chi)$ to 
$\C$ .
The only problem is to be assured that the residues of $L(s;\chi)$ are positive integers. 
But this follows from \eqref{eq:merologFormalSingularities} and the 
\textit{Gauss--Bonnet formula}
		\begin{equation*}
			\frac{\Vol(X)}{2\pi}= 2g-2,
		\end{equation*}
		where $g\geq2$ is the genus of the surface.
By \cite[Theorem 4.2.6]{frahm2023twisted}, we have the following result. 
	
	\begin{thm}\label{selbergmero}
		Let $X=\Gamma\backslash \H^2$ be a compact hyperbolic surface and let $\chi\colon \Gamma\rightarrow \GL(V_{\chi})$ be a finite-dimensional, complex representation of $\Gamma$. Then, the twisted Selberg zeta function $Z(s;\chi)$ admits a holomorphic continuation to $\C$ with zeros given by the following formal product:
		$$ \prod_{j=0}^\infty\left(s-\frac{1}{2}-i\mu_j\right)\left(s-\frac{1}{2}+i\mu_j\right)\prod_{k=0}^\infty(s+k)^{(2g-2)\dim(V_\chi)(1+2k)}. $$
	\end{thm}
	
The  meromorphic continuation to $\C$ follows from the \textit{relation} of the two 
dynamical zeta functions.  
By  \cite[Lemma 4.1.5]{frahm2023twisted}, we have 
\begin{equation*}
			R(s;\chi)=\frac{Z(s;\chi)}{Z(s+1;\chi)}.
		\end{equation*}
Hence, we get the following corollary (\cite[Corollary 4.2.8]{frahm2023twisted}).

	\begin{coro}\label{cor:RuelleMeromorphic}
		The twisted Ruelle zeta function $R(s;\chi)$ admits a meromorphic continuation to $\C$ with zeros and singularities given by the following formal product:
		$$ \prod_{j=0}^\infty\frac{(s-\frac{1}{2}-i\mu_j)(s-\frac{1}{2}+i\mu_j)}{(s+\frac{1}{2}-i\mu_j)(s+\frac{1}{2}+i\mu_j)}\times s^{(2g-2)\dim(V_\chi)}\prod_{k=1}^\infty(s+k)^{2(2g-2)\dim(V_\chi)}. $$
		In particular, $R(s;\chi)$ has a zero of order $(2g-2)\dim(V_\chi)$ at $s=0$.
	\end{coro}

The functional equation for the twisted Selberg zeta function follows also
from the resolvent trace formula  \eqref{finalresolvent}, 
by letting $s=s_1$ and considering the transform $s\mapsto 1-s$.
Then, since the left-hand side of \eqref{finalresolvent} remains invariant under this transform,
and using a classical identity for the term $\pi\tan\pi x $, we get the following
theorem (\cite[Theorem 4.2.9]{frahm2023twisted}).

	\begin{thm}
		The twisted Selberg zeta function satisfies the following functional equation.
		\begin{equation}\label{fefinal}
			\eta(s;\chi) = \frac{Z(s;\chi)}{Z(1-s;\chi)} = \exp\bigg[\dim(V_{\chi})\Vol(X)
			\int_0^{s-\frac{1}{2}}r\tan\pi r\,dr\bigg],
		\end{equation}
		where the integral is a complex line integral along any curve from $0$ to $s-\frac{1}{2}$.
	\end{thm}
	
	We remark that the value of the integral in \eqref{fefinal} depends on the chosen curve from $0$ to $s-\frac{1}{2}$, but the right-hand side of \eqref{fefinal} does not, because the residues of $r\tan\pi r$ are integer multiples of $\frac{1}{2\pi}$ and $\dim(V_\chi)\Vol(X)=2\pi(2g-2)\dim(V_\chi)$.

The following corollary (\cite[Corollary 4.2.10]{frahm2023twisted})
follows from the functional equation 	\eqref{fefinal} and  \cite[Lemma 4.1.5]{frahm2023twisted},

	\begin{coro}
		The twisted Ruelle zeta function satisfies the following functional equation
		\begin{equation*}
			R(s;\chi)R(-s;\chi)=(2\sin\pi s)^{2(2g-2)\dim V_\chi}.
		\end{equation*}
	\end{coro}

	The functional equation for $R(s;\chi)$ immediately gives the order of vanishing of $R(s;\chi)$ at $s=0$. We have the following corollary (\cite[Corollary 4.2.11]{frahm2023twisted}).
	
	\begin{coro}\label{Ruellezerosurf}
		The behaviour of the twisted Ruelle zeta function $R(s;\chi)$ near $s=0$
		is given by
		\begin{equation}\label{eq:VanishingOrderRuelle}
			R(s;\chi)=\pm(2\pi s)^{\dim(V_{\chi})(2g-2)}+\text{higher order terms}.
		\end{equation}
	\end{coro}

We remark here that the sign in \eqref{eq:VanishingOrderRuelle} can be determined explicitly.
\begin{equation*}
R(s;\chi) \sim (-1)^m(2\pi s)^{\dim(V_\chi)(2g-2)}+\text{higher order terms},
\end{equation*}
where $m$ denotes the multiplicity of the eigenvalue $\lambda=0$ of $\Delta_\chi^\sharp$.
For more details, we refer the reader to \cite[Remark 4.2.12]{frahm2023twisted}.

\subsection{Special values}

For a compact hyperbolic surface $X$, by Corollary \ref{Ruellezerosurf},
the order of vanishing of $R(s;\chi)$ at $s=0$,
is related to the \textit{Euler characteristic} 
$\chi(X)=2-2g$ of $X$.

For compact hyperbolic orbisurfaces, the behaviour
of $R(s;\rho)$ at $s=0$ is given by Theorem \ref{orbispecialvalue}.
Under certain assumptions for the representation $\rho$,
$R(s;\rho)$ at $s=0$ is equal to the \textit{Reidemeister--Turaev torsion} of 
the unit tangent bundle (up to a sign).
Again, we use the \textit{functional equations for the twisted 
Selberg zeta function} (see \cite[Section 4]{JEP_2023__10__1391_0}).

We mention here that in recent paper (\cite{jorgenson2026determinantstwistedlaplacianstwisted}),
Jorgenson, Smajlovic and the author proved, in the orbisurface case, a relation 
between the \textit{twisted Selberg zeta function}  $Z(s;\rho)$
associated with an arbitrary, finite-dimensional representation $\rho$,
and the \textit{regularized determinant of the twisted Laplacian} associated with $\rho$,
that involves a multiplicative constant, which is called the \textit{torsion factor}.

By \cite[Theorem 6.1]{jorgenson2026determinantstwistedlaplacianstwisted}, we have 
\begin{equation}
\mathrm{det}(\Delta_{\tau_m,\rho}^{\sharp} +s(s-1))= Z(s;\rho)Z_I(s,\rho)Z_{\rm ell}(s;\rho)e^{\widetilde C},
\end{equation}
where the functions $Z_{I}(s;\rho)$ and $Z_{\rm ell}(s;\rho)$ stem from the identity and elliptic contributions
to the \textit{trace formula}. The constant $\widetilde C$, the \emph{torsion factor}, is given by
\begin{equation}
\widetilde{C}=-\mathrm{dim}(V_\rho)\chi(X)(2\zeta'(-1)-\log\sqrt{2\pi})+ \sum_{j=1}^{r}\frac{\log\nu_j}{2\nu_j}\sum_{k=1}^{\nu_j-1}\frac{\mathrm{Tr}
  (\rho(\gamma_j)^k) e^{i\pi k m/\nu_j}}{\sin^2\left(\frac{\pi k}{\nu_j}\right)}
\end{equation}
where the sum is taken over representatives $\gamma_j$ of the $r$ inconjugate
elliptic classes with orders $\nu_j$, $\chi(X)$ is the Euler characteristic of $X$, and $\zeta$ denotes the Riemann zeta function (see \cite[Section 6]{jorgenson2026determinantstwistedlaplacianstwisted}).

In the case $\rho$ is the trivial character, the authors in 
\cite{jorgenson2026determinantstwistedlaplacianstwisted})
reproved by different means some of the main results of \cite{sarnak1987determinants}.

Moreover, by \cite[Corollary 8.1.1]{jorgenson2026determinantstwistedlaplacianstwisted}, 
we have that if $\mu=0$ is \textit{not} the eigenvalue of the twisted Laplacian, 
then the  \textit{twisted Selberg zeta function}
$Z(s;\rho)$ is non-vanishing at $s=1$ and
  \begin{equation*}
  \mathrm{det}\Delta_{\tau_m,\rho}^{\sharp}= Z(1;\rho)Z_I(1;\rho)Z_{\rm ell}(1;\rho)e^{\widetilde C}.
  \end{equation*}

In \cite{jorgenson2026determinantstwistedlaplacianstwisted} also, 
the authors studied the \textit{asymptotic behaviour of the torsion factor}
for a sequence of non-unitary representations introduced by Yamaguchi in \cite{Ya17} and prove that the asymptotic behaviour of this constant as the dimension of the representation tends to infinity is the same as the behaviour of the higher-dimensional Reidemeister torsion on the unit tangent bundle
(up to an absolute constant). In  particular, by \cite[Corollary 7.2]{jorgenson2026determinantstwistedlaplacianstwisted}, we have

\begin{coro}
Let $\{\rho_{N}\}_{N\geq 1}$ be a sequence of representations, such that $\mathrm{Tr}
  (\rho_{N}(\gamma_j))$ is uniformly bounded by some constant $C(X)$ which is independent of $j$ and $N$. Assume also that the sequence of dimensions $\{\dim (V_{\rho_N})\}_{N\geq 1}$ of representations $\rho_N$ is such that $\lim_{N\to\infty} \dim (V_{\rho_N})=\infty$. 
We then have that
  \begin{equation*}
\lim_{N\to\infty}\frac{ \widetilde C}{\dim (V_{\rho_N})}=-\chi(X)\left(2\zeta'(-1) - \log\sqrt{2\pi} \right).
\end{equation*}
  \end{coro}

In the case of compact hyperbolic odd-dimensional manifolds
and for specific representations, described in 
Subsection \ref{Friedodd},
the twisted Ruelle zeta function at zero is equal to 
a \textit{spectral invariant}, the \textit{Cappell-Miller torsion}.
In this case, the Cappell-Miller torsion
 "arise" from the \textit{spectral side} of the trace formula
(under certain technicalities).

\bibliographystyle{amsplain}
\bibliography{ref}

\vspace{5mm}\noindent
Polyxeni Spilioti\\
Department of Mathematics \\
University of Patras\\
Panepistimioupoli Patron 265 04\\
Greece\\
e-mail: pspilioti@upatras.gr

\end{document}